\documentclass[11pt]{amsart}

\usepackage[square,compress,comma, numbers,sort]{natbib}
\usepackage{amsfonts,mathtools,natbib}

\allowdisplaybreaks[4]

\usepackage{amsmath}
\usepackage{bbm}
\usepackage{amssymb}
\usepackage{mathtools}

\usepackage{color}

\def\zE#1{{\textcolor{c30}{#1}}}
\def\zE#1{#1}

\definecolor{c20}{rgb}{0.,0.7,0.}
\definecolor{c30}{rgb}{0.,0.,1.}
\definecolor{c40}{rgb}{1,0.1,0.7}
\definecolor{c50}{rgb}{1,0,0}
\definecolor{c60}{rgb}{1,0.9,0.1}

\def\ee#1{{\textcolor{cyan}{#1}}}
\def\ee#1{#1}

\def\cE#1{{\textcolor{cyan}{#1}}}
\def\cE#1{#1}

\def\kkk#1{{\textcolor{magenta}{#1}}}
\def\kkk#1{#1}

\newcommand{\ve}{\varepsilon}

\newcommand{\E}[1]{\mathbb{E}\left\{ #1\right\}}

\newcommand{\pk}[1]{\mathbb{P} \left\{ #1 \right \} }

\newcommand{\R}{\mathbb{R}}

\newcommand{\inr}{\in \R}

\newcommand{\limit}[1]{\lim_{#1 \to   \infty}}

\newcommand{\BQN}{\begin{eqnarray}}
\newcommand{\EQN}{\end{eqnarray}}
\newcommand{\BQNY}{\begin{eqnarray*}}
\newcommand{\EQNY}{\end{eqnarray*}}

\newcommand{\BS}{\begin{sat}}
\newcommand{\ES}{\end{sat}}
\newcommand{\BT}{\begin{theo}}
\newcommand{\ET}{\end{theo}}
\newcommand{\BK}{\begin{korr}}
\newcommand{\EK}{\end{korr}}
\newcommand{\EQD}{\stackrel{d}{=}}

\newcommand{\BD}{\begin{de}}
\newcommand{\ED}{\end{de}}

\newcommand{\BIT}{\begin{itemize}}
\newcommand{\EIT}{\end{itemize}}
\newcommand{\BDI}{\begin{description}}
\newcommand{\EDI}{\end{description}}

\newcommand{\BRM}{\begin{remarks}}
\newcommand{\ERM}{\end{remarks}}

\newcommand{\BEL}{\begin{lem}}
\newcommand{\EEL}{\end{lem}}

\newtheorem{theo}{Theorem}[section]
\newtheorem{sat}[theo]{Proposition}
\newtheorem{de}[theo]{Definition}
\newtheorem{lem}[theo]{Lemma}

\newtheorem{example}[theo]{Example}
\newtheorem{korr}[theo]{Corollary}
\newtheorem{remark}[theo]{Remark}
\newtheorem{remarks}[theo]{Remarks}

\newcommand{\nelem}[1]{{Lemma \ref{#1}}}

\newcommand{\netheo}[1]{{Theorem \ref{#1}}}

\newcommand{\prooftheo}[1]{ \textsc{\bf Proof of Theorem} \ref{#1}.}

\newcommand{\prooflem}[1]{\textsc{\bf Proof of Lemma} \ref{#1}.}

\newcommand{\COM}[1]{}

\def\td{\text{\rm d}}
\def\IF{\infty}

\newcommand{\QED}{\hfill $\Box$}

\newcommand{\kb}[1]{\boldsymbol{#1}}
\newcommand{\vk}[1]{\kb{#1}}

\def\ve{\varepsilon}

\def\IF{\infty}

\begin{document}

\title{Simultaneous Ruin Probability for Two-Dimensional Brownian and L\'evy Risk Models}

\author{Krzysztof D\c{e}bicki}
\address{Krzysztof D\c{e}bicki, Mathematical Institute, University of Wroc\l aw, pl. Grunwaldzki 2/4, 50-384 Wroc\l aw, Poland}
\email{Krzysztof.Debicki@math.uni.wroc.pl}

\author{Enkelejd  Hashorva}
\address{Enkelejd Hashorva, Department of Actuarial Science, 
University of Lausanne,\\
UNIL-Dorigny, 1015 Lausanne, Switzerland
}
\email{Enkelejd.Hashorva@unil.ch}

\author{Zbigniew Michna}
	\address{Zbigniew Michna, Department of Mathematics and Cybernetics, Wroc\l aw University of Economics, Poland}
	\email{zbigniew.michna@ue.wroc.pl}

\bigskip

\date{\today}
 \maketitle

 {\bf Abstract:} The ruin probability in the classical Brownian risk model can be explicitly calculated for both finite and
 infinite-time horizon. This is not the case for the simultaneous ruin probability in two-dimensional Brownian risk model. Resorting on asymptotic theory, we derive  in this contribution   approximations of both simultaneous ruin probability and simultaneous ruin time for the two-dimensional Brownian risk model  when the initial capital increases to infinity.
  Given the interest in proportional reinsurance, we consider in some details the case where the correlation is 1.
  This model is tractable allowing for explicit formulas for the simultaneous ruin probability for linearly dependent
  spectrally positive L\'evy processes. Examples include perturbed Brownian and gamma L\'evy processes.




\section{Two-Dimensional Brownian Risk Model}
The classical Brownian risk model (BRM) of an insurance portfolio
\[R_1(t)= u+ c_1t - \sigma_1W_1(t), \ \ \ \ \ t\ge 0,\]
with $W_1$ a standard Brownian motion, $\sigma_1>0$,
the initial capital $u>0$ and the premium rate $c_1>0$,
is a key benchmark model in risk theory; see e.g., \cite{ig:69}.

The ruin probability in the time horizon $[0,T]$ for some finite positive $T$ is given by (see e.g., \cite{DeM15})
\BQN \label{nuk}
\widetilde{ \psi}_T(u)\coloneqq\pk{ \inf_{t\in [0,T]} R_1(t) < 0}&=& \pk{\sup_{t\in [0,T]} (\sigma_1W_1(t)- c_1 t)> u}\notag \\
&=& \Phi\left(-\frac u{\sigma_1 \sqrt{T}} -\frac{c_1\sqrt{T}}{\sigma_1}\right)+
e^{-2c_1u/\sigma_1}\Phi\left(- \frac u{\sigma_1 \sqrt{T}} +\frac{c_1\sqrt{T}}{\sigma_1}\right)
\EQN
for any $u\geq 0$, with $\Phi$ the distribution function of an $N(0,1)$ random variable.\\
In the infinite-time horizon, i.e., for $T=\infty$, the corresponding ruin probability
for this risk model is
\BQN
\widetilde{ \psi}_\IF(u)\coloneqq\pk{ \inf_{t\ge 0} R_1(t) < 0}= e^{- 2c_1u/\sigma_1^2} .
\EQN
  Since in practice an insurance  company runs multiple portfolios simultaneously,
  it is of interest to calculate the simultaneous ruin probability for the classical benchmark BRM.
 For notational simplicity, we shall consider  only the two dimensional setup, where for the second  portfolio  we consider the risk process
$$R_2(t)= v+ c_2t- \sigma_2W_2(t), \quad t\ge 0,$$
with $W_2$ another standard Brownian motion, $v$ the initial capital and $c_2>0$.
Hereafter  $(W_1(t),W_2(t)),t\ge 0$ are assumed to be jointly Gaussian {with  the same law as
\BQN \label{BB}
(B_1(t), \rho B_1(t)+ \rho^* B_2(t)), \quad t\ge 0, \quad \rho^*= \sqrt{1- \rho^2}, \quad \rho \in (-1,1],
\EQN
where $B_1,B_2$ are two independent standard Brownian motions.}
{Thus the correlation between $W_1(t)$ and $W_2(t)$ is $\rho$ for $t>0$.}
The special case $\rho=1$ will be discussed separately in Section \ref{s.prop}. In this bivariate risk model, tractable expressions for the simultaneous ruin probability are not available for both finite and infinity-time horizon.
Here we are concerned with the study of the ruin probability
in finite-time, which from practical point of view  is more natural.

In the  2-dimensional BRM  the probability of simultaneous ruin
of both portfolios in  the time period [0,T] is given by
\BQNY
\pk{\exists t\in [0,T]:  R_1(t)< 0,  R_2(t)< 0}
 &=&
 \pk{\exists_{t\in [0,T]}: \sigma_1W_1(t)- c_1 t> u, \sigma_2W_2(t)- c_2 t> v},
 \EQNY
which by self-similarity (time-scaling property) of Brownian motion
reduces to
\[
\pk{\exists_{t\in [0,1]}: W_1(t)- \frac{c_1\sqrt{T}}{\sigma_1} t> \frac{u}{\sigma_1\sqrt{T}},
W_2(t)- \frac{c_2\sqrt{T}}{\sigma_2} t> \frac{v}{\sigma_2\sqrt{T}}}.
\]
Consequently,  in order to simplify the presentation,
we shall consider in the following $T=1$, $\sigma_1=\sigma_2=1$ and define for
any $u,v$ non-negative  the simultaneous ruin probability as
\[
 \psi(u,v)\coloneqq\pk{\exists_{t\in [0,1]}: W_1(t)- c_1 t> u, W_2(t)- c_2 t> v}.
\]
The main findings of this contribution concern the approximation of
\BQNY
 \psi(u,au)
 &=& \pk{\exists_{t\in [0,1]}: W_1(t)- c_1 t> u, W_2(t)- c_2 t> au}
 \EQNY
 as $u\to \IF$, for any given constant  $a\in (-\IF, 1]$. Note that there is no restriction to consider only $a\le 1$ and in our model it is possible to deal  also with $a,c_1,c_2$ being negative. This reflects the fact that depending on the correlations between two portfolios, the need for initial capital $u$ and $v$ can be different. \\
 Clearly, the simplest possible model is when $W_1$ and $W_2$ are independent. Even in this model, it is not possible to calculate $\psi(u,au)$ explicitly. Since for independent Gaussian processes the main tools of asymptotic theory of those processes are still available, the asymptotic behaviour of $\psi(u,au)$ can be established by modifying the classical approach (i.e., using Gordon inequality, see  \cite{Krzyspp}[Prop. 3.6],  instead of the well-known Slepian inequality, see e.g., \cite{LifBook,Pit96}).\\
 If  $W_1,W_2$ are jointly Gaussian and dependent, then the  calculation of the simultaneous ruin probability  is much more difficult to deal with, since  there is no substitute for Gordon inequality and the current methodology cannot cover the approximation of  extremes of vector-valued dependent risk processes, see also discussion in Section 2. In order to understand the asymptotic behaviour of the simultaneous ruin probability  as the initial capital $u$ tends to infinity,
 we present next a sharp bounds for $\psi(u,v)$, which  also give some insights on the asymptotic
 approximation of the simultaneous ruin probability when $u$ tends to infinity.\\
 First, observe that for any $u, c_1,c_2$ we have a simple upper bound
 $$ \psi(u,au) \le \min \left( \pk{\sup_{t\in [0,1]} (W_1(t)- c_1t)>u},
 \pk{\sup_{t\in [0,1]} (W_2(t)- c_2t)>au} \right)=:g(u,au).$$
 In view of \eqref{nuk} the  upper bound $g(u,au)$  can be calculated explicitly.
However, if $a\in (\rho,1]$ this upper bound is too rough as the next result shows. Throughout in the following
  $\mathbb{I}(\cdot)$ is the indicator function and $\rho^*= \sqrt{1- \rho^2} \in [0,1]$. Further $\Psi=1- \Phi$ with $\Phi$ the standard normal distribution on $\R$.

\BS \label{prop1}
For $c_1,c_2\in \R$, {$(u,v) \in \R^2\setminus (-\infty,0]^2$} and $ \rho \in (-1,1] $ we have
 \begin{eqnarray}
\label{ubb}
\pk{W_1(1)>u+c_1,W_2(1)>v+c_2}
\le
\psi(u,v)
\le
\frac{\pk{W_1(1)>u+c_1,W_2(1)>v+c_2}}{\pk{W_1(1)>\max(c_1,0),W_2(1)>\max(c_2,0)}}.
 \end{eqnarray}
 \ES

The main result of this contribution, given in \netheo{Th1} below, shows that a precise asymptotic approximation of $\psi(u,au)$, as $u\to \IF$,
can be obtained by using more advanced techniques.
\netheo{Th1} presents  interesting insight on the simultaneous probability of
ruin given the correlation $\rho$ that governs the risk processes $R_1$ and $R_2$.
 For this case, we have that if
the proportion of initial capitals between first and second risk \ee{process} is larger than the correlation function, \ee{that is $\rho \ge a$}, then the ruin probability is much smaller.

Related results for the infinite-time horizon are obtained in \cite{Mandjes07, PistAvram1, PistAvram2, Rolski17, MR3776537,borovkov17,
PalmBorov}.
The first three  papers consider the case that $\rho=1$. In \cite{Rolski17} the case $\rho \in (-1,1)$
is dealt with and \cite{borovkov17} extends \cite{PistAvram1,PistAvram2} to the $d$-dimensional setup of non-degenerated risk processes of Sparre-Andersen type.\\
 The asymptotic behaviour of the ruin probability in finite-time horizon,
 when $u\to \IF$, compared with the results of \cite{Rolski17} is completely different.
In particular, the leading term in the asymptotics  for the finite-time horizon
is $e^{- q_{a,\rho}u^2/2}$
with
 \BQN\label{qar}
q_{ a,\rho} = \frac{1\cE{-} 2a \rho + a^2}{1- \rho^2}\mathbb{I}(a>\rho)+ \mathbb{I}(a \le \rho).
\EQN
Note that if $a\in (\rho,1)$, then $q_{a,\rho}>1$.
In the infinite-time horizon, the leading term in the asymptotic of simultaneous ruin probability equals $e^{- c_{a,\rho} u}$ for some \ee{positive} $c_{a,\rho}$;
see \cite{Rolski17,PalmBorov}.

In the literature two-dimensional risk models are mainly concerned with heavy-tailed setup, see e.g.,
\cite{MR3102482,MR3493174, cl2017,foss2017two} and the references therein. The light-tailed assumption is different;
see \cite{AsmussenE} for some explanations and the difficulties in the light-tailed settings.

Brief  organisation of the rest of the paper. In the next section we give short discussions of our results including the case $\rho=1$ and the approximation of the conditional ruin time. All the proofs are displayed in Section 4.

\section{Main result}

 Let
  	   $\varphi_\rho$ stands for the joint probability density function (pdf) of $(W_1(1), W_2(1))$ and $\sim$  means asymptotic equivalence of two functions when the argument
  $u$ tends to infinity.
For
$a\in ( \rho,1]$ let the constant $C_{a,\rho}\in (0,\IF)$ be given by
\BQN
\quad  C_{a,\rho}= \int_{\R^2}
\pk{\exists_{t\ge 0 }:
	\begin{array}{ccc}
	W_1(t) - t > x \\
W_2(t) - at > y
	\end{array}
} e^{ \lambda_1 x+ \lambda_2 y }\, \td x \td y,
\label{car}
\EQN
where
\BQN \label{lam12}
\lambda_1=  \frac{1- a\rho}{1- \rho^2}, \quad \lambda_2=\frac{a- \rho }{1- \rho^2}
\EQN
are both positive.

\BT \label{Th1} Let $c_1,c_2$ be two given constants and let $\rho \in (-1,1)$.\\
 i) If $a\in (\rho,1]$, then  as $u\to \IF$
 \BQN
 \psi(u, au) \sim C_{a,\rho} u^{-2}\varphi_\rho(u+c_1,  a u+ c_2).
 \EQN
 ii) If $a\le \rho$, then we have as $u\to \IF$
 \BQN
 \psi(u, au) \sim  2\sqrt{2 \pi (1- \rho^2)} \Phi^*(c_1 \rho- c_2) e^{ \frac{ (c_2- \rho c_1)^2}{2(1- \rho^2)}} u^{-1}\varphi_\rho(u+c_1,  \rho u+ c_2)  ,
 \EQN
 where $\Phi^*(c_1 \rho- c_2)=1$ if $a< \rho$ and $\Phi^*$ is the df of $\sqrt{1- \rho^2} W_1(1)
 $  when $a=\rho$.\\
 iii) Note that $C_{a, \rho}$ is not a type of Pickands constant in these setting, see \cite{Rolski17} for the multivariate version of those constants and \cite{Pit96,MR3745392,KDEH1}.
  	\ET

\begin{remark} i) In view of \cite{MRMGaussH}{Theorem 4.1} {and Theorem \ref{Th1}}
we have that for $a\in (\rho ,1]$
	$$  \psi(u, au) \sim  C_{a,\rho} \lambda_1 \lambda _2 \pk{ W_1(1) > u+ c_1, W_2(1)>  a u+ c_2}, \quad u\to \IF,$$
	with
{ $\lambda_1,\lambda_2$ defined in (\ref{lam12}),}
whereas if $a \le \rho$,  then
$$  \psi(u, au) \sim   2 \pk{ W_1(1)  > u+ c_1, W_2(1)> a u+ c_2} \sim 2 \Phi^*(c_1 \rho- c_2) \pk{W_1(1)> u + c_1} , \quad u\to \IF. $$ 		
Moreover, {combination of Theorem \ref{Th1} with Proposition \ref{prop1}
gives
the following upper bound}
$$  C_{a,\rho} \le \frac{1}{ \lambda_1 \lambda _2 \pk{W_1(1)> \max(0,c_1), W_2(1)> \max(0,c_2)}}.$$
ii) From the above results, for any $a \in (\rho, 1]$ and $b \le \rho$ we have
$$ \limit{u} \frac{\psi(u, au)}{\psi(u,bu)}=0.$$
In particular, if $\rho=0$, the above holds for any $a\in (0,1], b\le 0$.
 \def\zE#1{{\textcolor{c30}{#1}}}  	
\end{remark}
Theorem \ref{Th1} enables us to analyze the simultaneous ruin time $\tau_{sim}(u)$ on $[0,1]$  defined by
$$ \tau_{sim}(u)= \inf\{ t\in [0,1]: W_1(t)- c_1 t> u, W_2(t)-c_2 t> au \}.$$
Our result below shows that  $
u^2(1-\tau_{sim}(u))
$ conditioned that $\tau_{sim}(u)\le 1$,
converges as $u\to\infty$, to an exponentially distributed random variable.
\BT\label{Th_3} If $a \le 1, \rho \in (-1,1) $ and $x\ge 0$, then with $q_{{a},\rho}$ defined in \eqref{qar} we have
\[
\lim_{u\to\infty}\pk{u^2(1-\tau_{sim}(u))\le x|\tau_{sim}(u)\le 1}=
1-\exp\left(-q_{{a}, \rho} x /2\right).
\]
\ET
Note that if $a> \rho$, then $q_{{a}, \rho} >1$ and $q_{{a},\rho}=1$ for $a\le \rho$.

\section{Proportional portfolios with one-sided L\'evy risk processes}\label{s.prop}
In this section we consider the  case when the insurance companies share the same portfolio of
claims, with some proportion $r_1,r_2>0$, respectively and {the portfolio is modeled by a L\'evy process}. This is typical for proportional reinsurance treaties.
We refer to, e.g., \cite{fumwe, Mandjes07, PistAvram1, PistAvram2} for the analysis of this model
for infinite-time ruin problem in Brownian and L\'evy setup.
Following recent results of Michna \cite{Mic17}, we shall derive exact distribution of
the corresponding ruin probability
for the claim process modeled by
a spectrally one-sided L\'evy process $Z$ with absolutely continuous one-dimensional distributions.
Since for \ee{positive} $r_1,r_2$ 
\[
 \pk{\exists_{t\in [0, T]}: r_1Z(t)- c_1 t> x, r_2Z(t)- c_2 t> y}
 =
 \pk{\exists_{t\in [0, T]}: Z(t)- \frac{c_1}{r_1} t> \frac{x}{r_1}, Z(t)- \frac{c_2}{r_2} t> \frac{y}{r_2}},
\]
in the rest of this section, with no loss of generality, we
suppose that $r_1=r_2=1$.
Thus the aim of this section is to obtain exact (non-asymptotic) expressions for the simultaneous ruin probability on finite time horizon $[0,T]$ defined {by}
$$
\psi_Z(x,y)= \pk{\exists_{t\in [0,T]}: Z(t)- c_1 t> x, Z(t)- c_2 t> y}.
$$
%
Below we exclude the degenerated scenario $c_1=c_2$
and by the symmetry of the considered problem we assume that $c_1>c_2$.
Utilising the findings in Michna \cite{Mic17} we shall derive an explicit formula for $\psi_Z(x,y)$
both for spectrally positive and spectrally negative $Z$,
which is the main result of this section. 

Suppose first that $Z$ is spectrally positive.
For $T,u$
positive and arbitrary constant $c$ set
\begin{eqnarray}\label{calL}
{\mathcal{L}}(c, T, u)\coloneqq\pk{Z(T)-cT>u}-\int_0^T\frac{\E{ \min(0,Z(T-s)-c(T-s))}}{T-s}\,
f(u+cs, s)\,{\rm d}s,\label{prob.L}
\end{eqnarray}
where $f(u,t)$ is the density function of $Z(t)$.
We note that in the light of  \cite{Mic17}, for $u\ge0$,
$$
\pk{\sup_{t\in [0,T]} (Z(t)-ct)>u}={\mathcal{L}}(c, T, u).
$$


\BT\label{Th_0}
Let
$Z$ be  a spectrally positive L\'evy process with c\'adl\'ag sample paths and $\pk{Z(0)=0}=1$.
Suppose that $Z(t),t>0$ has density function $f(u,t)$ and let  $c_1,c_2$ be {two given} constants such that
$\delta:=c_1- c_2> 0$.
\\
 i) If $x\geq y\ge 0$, then 
\BQN \label{rusl}
\hspace{ 1 cm} \mbox{}  \psi_Z(x,y)=
{\mathcal{L}}(c_1, T, x).
\EQN
ii) If $0 \le x<y <x+\delta T$, then
setting $\xi= (y-x)/\delta$
we have
\begin{eqnarray*}
\psi_Z(x,y)&=&
\mathcal{L}\left(c_2, \xi, y\right)+\int_{0}^\infty \mathcal{L}\left(c_1, T-\xi , z\right)
f(y+c_2\xi-z, \xi)\td z\\
&& -\int_{0}^\infty z \mathcal{L}\left(c_1, T- \xi , z\right)\td z  \int_{0}^{\xi}\frac{f(y+c_2 s,s)}{\xi-s}f\left(c_2\left(\xi-s\right)-z,\xi-s\right)\td s.
\end{eqnarray*}
iii) If $y\geq x+\delta T$ and $x\ge 0$, then \eqref{rusl} holds substituting $c_1, \ee{x}$ by $c_2,y$, respectively.
\ET

Next, let us suppose that the L\'evy process $Z$  is spectrally negative.
In view of \cite{Mic17}[Thm 5] we obtain  the following result.
\BT\label{Th_01}
Let
$Z$ be  a spectrally negative L\'evy process with c\'adl\'ag sample paths and $\pk{Z(0)=0}=1$.
Suppose that $Z(t),t>0$ has density function $p(u,t)$ and let  $c_1,c_2$ be {two given} constants such that
$\delta:=c_1- c_2> 0$.
\\
i) If $x\geq y\ge 0$, then 
\BQN \label{rus2}
\hspace{ 1 cm} \mbox{}  \psi_Z(x,y)=
x\int_0^T\frac{p(x+c_1s,s)}{s}\,\td s.
\EQN
ii) If $0 \le x<y <x+\delta T$, then
setting $\xi= (y-x)/\delta$
we have
\begin{eqnarray*}
\psi_Z(x,y)&=&
y\int_0^{\xi}\frac{p(y+c_2s,s)}{s}\,\td s
+
\int_{0}^\infty
z p(-z+y+c_2 \xi,\xi)\td z \int_0^{T-\xi}\frac{p(z+c_1s,s)}{s}\,\td s \,
\\
&&-y\,
\int_{0}^\infty
z \td z\int_0^{T-\xi}\frac{p(z+c_1s,s)}{s}\,\td s
\\
&&\,\,\,\,\,\,\cdot\int_{0}^{\xi}\frac{p(-z+c_2t, t)}{\xi-s}\,p(u+c_2(\xi-t), \xi-t)\td t\,
\end{eqnarray*}
iii) If $y\geq x+\delta T$ and $x\ge 0$, then \eqref{rusl} holds substituting
$c_1, \ee{x}$ by $c_2,y$, respectively.
\ET

In the rest of this section we apply Theorem \ref{Th_0}
to {important  L\'evy risk models}.
\begin{example}
If $Z(t),t\ge 0$ is a standard Brownian motion, then Theorem \ref{Th_0} is satisfied with
$f(u,t)=\frac{1}{\sqrt{2\pi t}}e^{-\frac{u^2}{2t}}$ and
\begin{equation}\label{ab}
\mathcal{L}(c,T, u)=\Phi(-uT^{-1/2}-c\sqrt{T})+e^{-2uc}\Phi(-uT^{-1/2}+c\sqrt{T}).\nonumber
\end{equation}
\end{example}
\begin{example}\label{ex2}
Let $Z$ be a  {\it gamma L\'evy} process with parameter $\lambda>0$ where the density function of $Z(t), t>0$  is given by
$$f(u, t)=\frac{\lambda^t}{\Gamma(t)}u^{t-1}e^{-\lambda u}\mathbb{I}_{\{u>0\}}.$$
Then Theorem \ref{Th_0} holds with
\begin{eqnarray*}
\mathcal{L}(c, T, u)&=&
\frac{\lambda^T}{\Gamma(T)}\int_{u+cT}^\infty z^{T-1}e^{-\lambda z}\,\td z\\
&&+
\,\lambda^T e^{-\lambda u} \int_0^T \td s\int_0^{c(T-s)}\frac{(u+cs)^{s-1}e^{-c\lambda s }}{\Gamma(s)\Gamma(T-s+1)} (c(T-s)-z)z^{T-s-1}e^{-\lambda z} \td z
\end{eqnarray*}
for $c,T,u$ positive.
\end{example}
\begin{example}
Suppose that $Z=Z_{\alpha,1,1}$ is an {\it $\alpha$-stable L\'evy} process with
$1<\alpha<2$, $\beta=1$ (i.e., skewed to the right) and scale parameter $\sigma=1$;
see, e.g.,  Samorodnitsky and Taqqu \cite{SaT}.
Then
$$f(u,t)=
\frac{1}{\pi t^{1/\alpha}}
\int_0^\infty e^{-x^\alpha}\cos\left(uxt^{-1/\alpha}-x^\alpha\tan{\frac{\pi\alpha}{2}}\right)\td x
$$
and Theorem \ref{Th_0} is satisfied with
\begin{eqnarray*}
\mathcal{L}(c,T,u)&=&
\frac{1}{\pi T^{1/\alpha}} \int_u^\infty\td z
\int_0^\infty e^{-x^\alpha}\cos\left((z+cT)xT^{-1/\alpha}-x^\alpha\tan{\frac{\pi\alpha}{2}}\right)\td x\nonumber\\
&& -\frac{1}{\pi}\int_0^T \frac{\E{\min(0,Z_{\alpha,1,1}(T-s)-c(T-s)) }}{(T-s)s^{1/\alpha}} \td s  \int_0^\infty e^{-x^\alpha}\cos\left((u+cs)xs^{-1/\alpha}-x^\alpha\tan{\frac{\pi\alpha}{2}}\right) \td x
\end{eqnarray*}
for $T>0$, $c\in \R$ and $u>0$,
where
$$
\E{\min(0,Z_{\alpha,1,1}(s)-cs)}=\frac{1}{\pi s^{1/\alpha}}  \int_{-\infty}^0 z\td z\int_0^\infty e^{-x^\alpha}\cos\left((z+cs)xs^{-1/\alpha}-x^\alpha\tan{\frac{\pi\alpha}{2}}\right) \td x.
$$
 \end{example}
\begin{example}
Consider gamma L\'evy risk process perturbed by Brownian motion, i.e.
suppose that
$Z(t)=Z_1(t)+\sigma Z_2(t)$, where $Z_1(t),t\ge0$ is a gamma L\'evy process, as defined in Example \ref{ex2},
$Z_2(t),t\ge 0$ is a standard Brownian motion independent of $Z_1$ and $\sigma>0$.
Then Theorem \ref{Th_0} holds with
\[
f(u,t)=\frac{\lambda^t}{\Gamma(t)\sigma\sqrt{2\pi t}}\int_0^\infty e^{-\frac{(x-y)^2}{2\sigma^2 t}-\lambda y}y^{t-1}\td y
\]
and
\begin{eqnarray*}
\mathcal{L}(c,T,u)& =&
\frac{\lambda^T}{\Gamma(T)\sigma\sqrt{2\pi T}}\int_u^\infty\td z
\int_0^\infty e^{-\frac{(z+cT-y)^2}{2\sigma^2 T}-\lambda y}y^{T-1}\td y \nonumber\\
&& -\frac{1}{\sigma\sqrt{2\pi}}\int_0^T \frac{\E{\min(0,Z(T-s)-c(T-s))} }{(T-s)\Gamma(s)\sqrt{s}}\lambda^s \td s  \int_0^\infty e^{-\frac{(u+cs-y)^2}{2\sigma^2 s}-\lambda y}y^{s-1}\td y
\end{eqnarray*}
for $T>0$, $c\in \R$ and $u>0$,
where
$$
\E{\min(0,Z(s)-cs)}=\frac{\lambda^s}{\Gamma(s)\sigma\sqrt{2\pi s}}\int_{-\infty}^0 z\td z\int_0^\infty e^{-\frac{(z+cs-y)^2}{2\sigma^2 s}-\lambda y}y^{s-1}\td y\,.
$$
\end{example}

\section{Proofs}
First recall that in our notation $B_1,B_2$ are two independent standard Brownian motions
and  $(W_1,W_2)$ has law given by \eqref{BB} for some $\rho \in (-1,1)$.
In order to shorten the notation, in the following we set
$W_i^*(t)= W_i(t)- c_it, i=1,2$, with  $c_1,c_2$ two given constants (not necessarily positive).
We shall write $\Psi_\rho$ for the tail distribution function of  $(W_1(1), W_2(1))$ and
 $\varphi_\rho$ for its pdf.

\subsection{Proof of Proposition \ref{prop1}}
The proof of the lower bound is immediate.
For the proof of the upper bound we follow the same idea as in the proof of
\cite{DDEL17}[Thm  1.1].
We shall use the standard notation for vectors which are denoted in bold.
Let $\vk{W}(t)=(W_1(t), W_2(t))$, $\vk{c}=(c_1,c_2)$.
For $\vk{u}=(u_1,u_2){\in \R^2\setminus(-\infty,0]^2}$ define next
$\vk{B}_{\vk{u}}:= \{(x,y)\in \R^2: (x=u_1 \wedge y\ge u_2)\vee (x\ge u_1 \wedge y=u_2)\}$,
{the boundary of the set $\{(x,y)\in \R^2: x\ge u_1,  y\ge u_2\}$,}
and
$$\tau(\vk{u}):=\inf\{t {\in [0,1]}:W_1(t)-c_1t\ge u_1, W_2(t)-c_2t\ge u_2  \}.$$
Observe that
\begin{eqnarray*}
\lefteqn{
\pk{W_1(1)-c_1\ge u_1, W_2(1)-c_2\ge u_2}}\\
&=&
\int_0^1
\pk{\tau(\vk u)\in \td t}
\int_{\vk{B}_{\vk{u}}}
\pk{\vk{W}(t)-\vk{c}t\in \td\vk{x}|\tau(u)=t}
\pk{\vk{W}(1-t)-\vk{c}(1-t)\ge \vk{u}-\vk{x}},
\end{eqnarray*}
where we used the Strong Markov property
{and the fact that $\vk W(1)- \vk W(t)$ has the same law as $\vk W(1-t)$ for any $t\in [0,1]$.}
Since for $x_1\ge  u_1, x_2 \ge u_2$ and any $t\in [0,1]$ we have
\begin{eqnarray*}
\pk{\vk{W}(1-t)-\vk{c}(1-t)\ge \vk{u}-\vk{x}}
&\ge&
\pk{\vk{W}(1-t)\ge\vk{\tilde{c}}(1-t)}\\
&=&
\pk{\vk{W}(1)\ge\vk{\tilde{c}}\sqrt{1-t}}\\
&\ge&
\pk{\vk{W}(1)\ge\vk{\tilde{c}}},
\end{eqnarray*}
where
$\vk{\tilde{c}}=(\max(c_1,0),\max(c_2,0))$, the proof is complete.
\QED

\subsection{Proof of Theorem \ref{Th1}}

Let in the following $\delta(u,T)=1-Tu^{-2}$, for $T,u>0$.
Before proceeding to the proof of Theorem \ref{Th1}
we present two lemmas: Lemma \ref{l.negl} that
provides a sharp upper bound for
{
\[
m(u,T)\coloneqq\pk{\exists_{t \in [0,\delta(u,T)]}: W_1(t) - c_1 t> u, W_2(t)- c_2 t> au }
\]
}
and Lemma \ref{singA} which gives
precise asymptotics of
\[
M(u,T)\coloneqq\pk{\exists_{t \in [\delta(u,T),1]}: W_1(t)- c_1 t> u, W_2(t)- c_2 t> au }
\]
as $u \to \IF$.

\BEL\label{l.negl}
For any $T>0$, $a \in (-\IF, 1] $ and sufficiently large $u$
\begin{eqnarray*}
	{m(u,T)} &\le &
e^{-T/{8}}
\frac{\pk{W_1(1)\ge u+ c_1, W_2(1)\ge au+ c_2}}{\pk{W_1(1)>\max(c_1,0),W_2(1)>\max(c_2,0)}}.
\end{eqnarray*}
\EEL

\prooflem{l.negl} For notation simplicity, we suppress the argument $u$ writing only $\delta(T)$ instead of $\delta(u,T)$ in the following.
By the self-similarity of Brownian motion, combined with Proposition \ref{prop1}, for any $u>0$
\begin{eqnarray*}
m(u,T)&{=}&
\pk{\exists_{t \in [0,1]}: W_1(t)-c_1\delta^{1/2}(T)t> \delta^{-1/2}(T)u,
                           W_2(t)-c_2\delta^{1/2}(T)t> \delta^{-1/2}(T)au }\nonumber\\
&\le&
\frac{\pk{W_1( \delta(T))\ge u+ c_1\delta (T), W_2(\delta(T)) \ge au +c_2 \delta(T) }}
     {\pk{W_1(1)>\max(c_1,0),W_2(1)>\max(c_2,0)}}.\label{eq11}
\end{eqnarray*}

\def\pp{\nu}
\def\ppy{\overline{\nu}}
Since, for sufficiently large $u$ (set below $\pp=\delta^{-1/2}(T), \ppy=\delta^{-1/2}(T/2)$ and recall that both $\pp$ and  $\ppy$ depend on $u$)
\begin{eqnarray*}
\pp u+ c_1/\pp
&\ge&
\ppy(u+ c_1), \quad
\pp au+ c_2/\pp
\ge
\ppy(au+ c_2),
\end{eqnarray*}
then we have
\begin{eqnarray*}
\lefteqn{
\pk{W_1(1)\ge \pp u+ c_1/\pp , W_2(1)\ge  \pp au +c_2/\pp  }}\\
&\le&
\pk{W_1(1)\ge \ppy (u+ c_1), W_2(1)\ge \ppy (au +c_2)}\\
&=&
{\ppy}^2
\int_{u+ c_1}^{\infty}
\int_{au+ c_2}^{\infty}
\varphi_\rho\left(\ppy x,\ppy y\right) \td x \td y
\end{eqnarray*}
 for sufficiently large $u$.
 Taking into account that
\begin{eqnarray*}
\varphi_\rho\left(\ppy x,\ppy y\right)
&=&
\frac{1}{2\pi (1-\rho^2)}
\exp\left( -\frac{{\ppy}^2 }{2(1-\rho^2)}(x^2-2\rho xy+y^2) \right)\\
&\le&
\frac{1}{2\pi (1-\rho^2)}
\exp\left( -\frac{1+Tu^{-2}{/2}}{2(1-\rho^2)}((\rho x- y)^2 +(1-\rho^2)x^2) \right)\\
&\le&
\varphi_\rho(x,y)
\exp\left( -\frac{Tu^{-2}}{{4}}x^2 \right)
\end{eqnarray*}
we get, for sufficiently large $u$ that
\begin{eqnarray*}
{\ppy}^2
\int_{u+ c_1}^{\infty}
\int_{au+ c_2}^{\infty}
\varphi_\rho\left(\ppy x,\ppy y\right) \td x \td y
&\le&
{\ppy}^2\int_{u+ c_1}^{\infty}
\int_{au+ c_2}^{\infty}
\varphi_\rho\left(x,y\right)
\exp\left( -\frac{Tu^{-2}}{{4}}x^2 \right)
\td x \td y
\\
&\le&
{ \exp\left( -\frac{T}{{8}} \right)}
\pk{W_1(1)\ge u+c_1, W_2(1)\ge au+c_2}.
\end{eqnarray*}
This completes the proof.
\QED

\BEL
i) For any $a\in (\rho, 1]$ and any $T>0$ we have 
\label{singA}
\BQN
M(u,T)&\sim & u^{-2} \varphi_\rho(u+c_1,au+c_2) I(T) , \quad u\to \IF,
\EQN
where
$$ I(T)\coloneqq \int_{\R^2}
\pk{\exists_{t\in [0, T ]}:
	\begin{array}{ccc}
	W_1(t) - t > x \\
	W_2(t) - at > y
	\end{array}
} e^{ \lambda_1  x + \lambda_2 y }\, \td x \td y \in (0,\IF).$$
ii) For any $a \le \rho,T>0$ with $\rho \in (-1,1)$ we have
\BQNY
M(u,T) \sim  u^{-1}  \varphi_\rho(u+ c_1,\rho u+ c_2)  I(T) ,  \quad u\to \IF,
\EQNY
where
$$I(T)\coloneqq\int_{\R^2}\pk{\sup_{t\in [0,T]} (W_1(t)-t)> x} \Bigl[\mathbb{I}(a <\rho)+ \mathbb{I}(y\zE{<} 0,a=\rho)\Bigr] e^{  x - \frac{y^2-2{y(c_2- c_1 \rho)}}{2 (1- \rho^2)}}\, \td x \td y. $$
\EEL

\prooflem{singA}
 For any $x,y$ and
$$u_x=u + c_1 - x/u, \quad u_y=au + c_2- y/u$$
we have (recall that  the pdf of $(W_1(1),W_2(1))$ is denoted $\varphi_\rho $)

$$\varphi_\rho(u_x , u_y ) =\colon
\varphi_\rho(u+c_1, au+c_2) \psi_u(x,y) \sim  \varphi_\rho(u+c_1, au+c_2)  e^{ \lambda_1 x + \lambda_2 y }, \quad u\to \IF,	 $$
where $\lambda_i$'s are given in \eqref{lam12}.
Set  below
$$ u_{x,y}:= u_y - \rho u_x= (a- \rho) u- (y-\rho x)/u+ c_2- {\rho}c_1$$
and let {$B_{1},B_{2}$} be two independent standard Brownian motions.

\COM{

The above convergence holds for almost all $x,y \inr$, consequently  we have
}

For any $u>0$ set further   $\bar t_u= 1- t/u^2$ and
\BQN \label{xu}
\bar x_u= 1- x/u^2, \quad x_u = x/u^2, \quad A(u)=  u^{-2} \varphi_\rho(u+c_1,au+c_2) .
\EQN
 $i)$  We have (recall $W_i^*(t)= W_i(t)- c_it$)
\BQNY
M(u,T)	&=&
u^{-2}	\int_{\R^2}
\pk{\exists_{t\in [\delta(u,T) ,1]}: W_1^*(t)>u, W_2^*(t)> au  \Bigl \lvert W_1(1)=u_x, W_2(1)= u_y  } \varphi_\rho(u_x,u_y)\, \td x \td y\\	
&=& A(u)
\int_{\R^2}
\pk{\exists_{t\in [0, T ]}:
	\begin{array}{ccc}
		B_1(\bar t_u )- c_1\bar t_u > u\\
		\rho B_1( \bar t_u ) + \rho^* B_2(\bar t_u )-  c_2 \bar  t_u > au
	\end{array}
	\Biggl \lvert
	\begin{array}{ccc} B_1(1)=u_x \\
		\rho^* B_2(1)= u_{x,y}
	\end{array}
} \psi_u(x,y)\, \td x \td y\\
&=\colon&A(u)
\int_{\R^2}  h_u(T,x,y) \psi_u(x,y)\, \td x \td y,
\EQNY
where {$\rho^*=\sqrt{1-\rho^2}$}. {For notational simplicity we define further 
$$ B_{u,1}(t) + \bar t_u u_x  \coloneqq B_{{1}}(\bar t_u  ) \lvert (B_{{1}}(1)= u_x) ,   \quad t\ge 0,$$
$$ B_{u,2}(t) +  \bar t_u   u_{{x},y}/\rho^*\coloneqq  B_2(\bar t_u ) \lvert ( \rho^* B_2(1)= u_{x,y}),   \quad t\ge 0.$$
}
The following weak convergence holds for all $x\inr $
$$u \Bigl[B_{u,1}(t) +  \bar t_u  u_x - c_1 \bar t_u  -u\Bigr] \to
B_1(t) - t
- x   , \quad t\in [0,T]$$
as $u\to \infty$. The above implies the weak convergence  as $u\to \IF${
$$u {\rho^*}  B_{u,2}(t) + u\Bigl[\bar t_u  u_{x,y} - (c_2-{\rho}c_1)\bar t_u  -u(a- \rho) \Bigr] \to
{\rho^*}B_2(t) - (a-\rho)t- (y-\rho x)   , \quad t\in [0,T]$$
for any $x,y\inr$. }
The following function
 $$ h(T,x,y)= \pk{ \sup_{t\in [0,T]} \min\Bigl( B_1(t)- t- x, B_{12}(t)- y\Bigr) > 0}$$
 is {non-increasing}
in both $x$ and $y$ and therefore it is continuous for $x,y\inr$ almost everywhere where
\BQN
\label{b12}
B_{12}(t)= 		\rho [B_{1}(t) - t ]  + \rho^*[
B_2(t) - t (a- \rho)/\rho^*].
\EQN
{
Note that by the independence of $B_1$ and $B_2$ we have that
$(B_1(t){-t}, B_{12}(t)), t\ge 0$  has the same law as
$$(W_1(t)- \textbf{}t,  W_2(t)- at),\quad t\ge 0$$
implying {that}
 $$ h(T,x,y)= \pk{ \sup_{t\in [0,T]} \min\Bigl( W_1(t)- t- x, W_{2}(t)- at- y\Bigr) > 0}>0.$$
}For any $(x,y)$ continuity point of $h(T,x,y)$, since $B_1,B_2$ are independent it follows by continuous mapping theorem that
\BQNY
h_u(T,x,y)&\coloneqq&\pk{\exists_{t\in [0, T ]}:
	\begin{array}{ccc}
		B_{u,1}(t) + \bar t_u  u_x  - c_1\bar t_u > u\\
		\rho [B_{u,1}(t) + \bar t_u u_x ]  + \rho^*[
		B_{u,2}(t) + \bar t_u u_{x,y}/\rho^*]  -c_2\bar t_u  > au
	\end{array}}\\
&=& \mathbb{P} \Bigl\{ \sup_{t\in [0, T ]} \min \Bigl(
		u(B_{u,1}(t) + \bar t_u  (u_x  - c_1) - u),
		u( \rho [B_{u,1}(t)  + \bar t_u u_x ] \\
		&&   + \rho^*[
		B_{u,2}(t) + \bar t_u u_{x,y}/\rho^*]  -c_2\bar t_u  - au)
		\Bigr)> 0 \Bigr\}\\
&\to & h(T,x,y), \quad u\to \IF\,.
\EQNY

The above convergence holds for almost all $x,y \inr$, consequently  using the dominated convergence theorem, we have
\BQNY
M(u,T)	&=&A(u)
\int_{\R^2}  h_u(T,x,y) \psi_u(x,y)\, \td x \td y\\
&\sim &A(u) \int_{\R^2} h(T,x,y) e^{\lambda_1 x+ \lambda_2 y}\, \td x \td y, \quad u\to \IF.
\EQNY
The application of the dominated convergence theorem can be justified as follows.
First note that for all $u$ large and some $\ve>0$ we have
$$ \psi_u(x,y)\le e^{ \lambda_{1, \ve} x+ \lambda_{2, \ve} y} , \quad x,y \inr,$$
where $\lambda_{i,\ve}{(x)}= \lambda_i + sign(x) \ve$.
\kkk{Moreover, using that for sufficiently large $u$ and $s,t\in[0,T]$ we have
$u^2\E{\left(B_{u,i}(t)-B_{u,i}({s})\right)^2} \le {\rm Const}|t-s|$
 for some ${\rm Const}>0$,
the application of Piterbarg inequality (see, e.g., \cite{Pit96}[Thm 8.1]) implies that
for $x,y\ge0$ and some constant $C_1$
\begin{eqnarray*}
h_u(T,x,y)
&=&
\mathbb{P} \Bigl\{ \sup_{t\in [0, T ]} \min \Bigl(
		u(B_{u,1}(t) + \bar t_u  (u_x  - c_1) - u),
		u( \rho [B_{u,1}(t)  + \bar t_u u_x ]\\
		 &&  + \rho^*[
		B_{u,2}(t) + \bar t_u u_{x,y}/\rho^*]  -c_2\bar t_u  - au)
		\Bigr)> 0 \Bigr\}\\
&\le&
\mathbb{P} \Bigl\{ \sup_{t\in [0, T ]} \Bigl(
		u(B_{u,1}(t) + \bar t_u  (u_x  - c_1) - u)
+
		u( \rho [B_{u,1}(t)  + \bar t_u u_x ]\\
		&&   + \rho^*[
		B_{u,2}(t) + \bar t_u u_{x,y}/\rho^*]  -c_2\bar t_u  - au)
		\Bigr)> 0 \Bigr\}\\
&\le&
\zE{\mathbb{P} \Bigl\{ \sup_{t\in [0, T ]} \Bigl(
		uB_{u,1}(t) +u\rho B_{u,1}(t)  + u\rho^*B_{u,2}(t) - (t(a+1)+x+y-\epsilon)
		\Bigr)> 0 \Bigr\}}\\
&\le&
\zE{\mathbb{P} \Bigl\{ \sup_{t\in [0, T ]} \Bigl(
		uB_{u,1}(t) +u\rho B_{u,1}(t)  + u\rho^*B_{u,2}(t)
		\Bigr)> x+y-C_1 \Bigr\}}\\
&\le&
\bar{C}e^{- C (x+y)^2} \le \bar{C}e^{-C(x^2+y^2)},
\end{eqnarray*}
for $x\ge0,y\le0$
\begin{eqnarray*}
h_u(T,x,y)
\le
\pk{\sup_{t\in [0,T] } u(B_{u,1}(t) + \bar t_u  (u_x  - c_1) - u)>0 } \le \bar{C}e^{- C x^2}
\end{eqnarray*}
and for $x\le0,y\ge0$
\begin{eqnarray*}
h_u(T,x,y)
\le
\pk{\sup_{t\in [0,T] } u( \rho [B_{u,1}(t)  + \bar t_u u_x ]
		                  + \rho^*[B_{u,2}(t) + \bar t_u u_{x,y}/\rho^*]  -c_2\bar t_u  - au)> 0 } \le \bar{C}e^{- C y^2}
\end{eqnarray*}
for some $C,\bar{C}>0$.
Hence we have
\begin{eqnarray*}
\int_{ x\ge 0, y \ge 0}
h_u(T,x,y){\psi_u(x,y)}  \td x \td y
&\le&
\bar{C}{\int_{0}^\infty e^{ (\lambda_1+ \ve) x}e^{- C x^2}\td x  \int_0^\IF
{e^{ (\lambda_2+ \ve) y}}e^{- C y^2} \, \td y}<\infty\\
\int_{ x\ge 0, y \le 0}
h_u(T,x,y){\psi_u(x,y)}  \td x \td y
&\le&
\bar{C}{\int_{-\IF}^0 e^{ (\lambda_2- \ve) y}\td y  \int_0^\IF
{e^{ (\lambda_1+ \ve) x}}e^{- C x^2} \, \td x}<\infty\\
\int_{ x\le 0, y \ge 0}
h_u(T,x,y){\psi_u(x,y)}  \td x \td y
&\le&
\bar{C}{\int_{-\IF}^0 e^{ (\lambda_1- \ve) x}\td x  \int_0^\IF
{e^{ (\lambda_2+ \ve) y}}e^{- C y^2} \, \td y}<\infty\\
\int_{ x\le 0, y \le 0}
h_u(T,x,y){\psi_u(x,y)}  \td x \td y
&\le&
{\int_{-\IF}^0 e^{ (\lambda_1- \ve) x}\td x
\int_{-\IF}^0
e^{ (\lambda_2- \ve) y} \, \td y}<\infty,
\end{eqnarray*}
which confirms the validity of the dominated convergence theorem.}

$ii)$ Next, when $a \le \rho$ we shall apply a different transformation, namely
$$ u_x= u+ c_1- x/u, \quad u_y=\rho  u + c_2- y. $$
With this notation we have
$$ u_{x,y}\coloneqq u_y - \rho u_x=  \rho x/u-y+ c_2- {\rho}c_1$$
and
$$\varphi_\rho(u_x , u_y )=\colon \psi_u(x,y) e^{  x - \frac{y^2}{2 (1- \rho^2)}} \sim \varphi_\rho(u + c_1, \rho u+ c_2)
e^{  x - \frac{y^2 - 2{y(c_2- c_1 \rho)}}{2 (1- \rho^2)}}, \quad u\to \IF.	 $$
For any $x,y\inr $ we have  thus
\begin{eqnarray*}
\limit{u} h_u(T,x,y)
&=&\pk{\sup_{t\in [0,T]} (W_1(t)-t)> x} \Bigl[  \mathbb{I}(a <\rho)+ \mathbb{I}({y< 0},a=\rho)\Bigr]\\
&&{+\pk{ \sup_{t\in [0,T]} \min\Bigl( W_1(t)- t- x, W_{2}(t)- at\Bigr) > 0}\mathbb{I}(y=0, a=\rho)}\\
&
=\colon&h(T,x,y).
\end{eqnarray*}
Setting $A(u): =u^{-1} \varphi_\rho(u+ c_1,\rho u+ c_2)$  we have further
\BQNY
\lefteqn{M(u,T)=}\\
&& u^{-1}
\int_{\R^2}
\pk{\exists_{t\in [\delta(u,T),1]}: W_1^*(t)>u, W_2^*(t)> au  \Bigl \lvert W_1(1)=u_x, W_2(1)= u_y  } \varphi_\rho(u_x,u_y)\, \td x \td y\\	
&=&
{u^{-1}}\int_{\R^2}
\pk{\exists_{t\in [0, T ]}:
	\begin{array}{ccc}
		B_1(\bar t_u )- c_1\bar t_u > u\\
		\rho B_1( \bar t_u ) + \rho^* B_2(\bar t_u )-  c_2 \bar  t_u > au
	\end{array}
	\Biggl \lvert
	\begin{array}{ccc} B_1(1)=u_x \\
		\rho^* B_2(1)= u_{x,y}
	\end{array}
} \psi_u(x,y) {e^{  x - \frac{y^2}{2 (1- \rho^2)}}}\, \td x \td y\\
&\sim & A(u) \int_{\R^2} h(T,x,y)e^{  x - \frac{y^2-2{y(c_2- c_1 \rho)}}{2 (1- \rho^2)}}\, \td x \td y, \quad u\to \IF.
\EQNY
The application of the dominated convergence theorem is simpler in this case and is therefore omitted.
 \QED

\prooftheo{Th1} Recall first that we define $\delta(u,T)=1-Tu^{-2}$. In view of \nelem{l.negl}, combined with
  Proposition \ref{prop1}, we
immediately obtain that
\[
\lim_{T\to\infty}\lim_{u\to\infty} \frac{\pk{ \exists_{t\in [0,\delta(u,T)]}: W_1^*(t)> u, W_2^*(t)> au}}{\psi(u,au)} =0.
\]
Hence, using that (recall $M(u,T)\coloneqq\pk{ \exists_{t\in [\delta(u,T), 1]}: W_1^*(t)> u, W_2^*(t)> au}$)
\begin{eqnarray*}
M(u,T)\le \psi(u,au) \le
\pk{ \exists_{t\in [0,\delta(u,T)]}: W_1^*(t)> u, W_2^*(t)> au}
+M(u,T)
\end{eqnarray*}
we obtain
\[
\lim_{T\to\infty}\lim_{u\to\infty}
\frac{M(u,T)}{\psi(u,au)} =1.
\]
Consequently, in view of
Lemma \ref{singA},
it suffices to prove that
\[
\lim_{T\to\infty}I(T)\in (0,\infty),
\]
where $I(T)$ is defined in Lemma \ref{singA}.
We derive the above one considering separately $a\in(\rho,1]$ and $a\le\rho$.\\

i) If  $a\in (\rho, 1]$, then we have
\BQNY \lim_{T\to\infty}
I(T)&=&\limit{T} \int_{\R^2} h(T,x,y)e^{\lambda_1 x+ \lambda_2 y}\, \td x \td y \\
& =&
\int_{\R^2} \limit{T} h(T,x,y)e^{\lambda_1 x+ \lambda_2 y} \, \td x \td y \in (0,\IF),
\EQNY
where $h$ is as in the proof of  \nelem{singA},  $\lambda_1,\lambda_2$ are positive constants defined in \eqref{lam12} and
$$ \widetilde h(x,y)\coloneqq\limit{T} h(T,x,y).$$
We have the following upper bound   
\BQNY
\lefteqn{
  \int_{ \R^2} \widetilde h(x,y)  e^{\lambda_1 x+ \lambda_2 y}\, \td x \td y}\\
 &\le&  \sum_{i=0}^{\IF } \int_{\R^2} h([i,i+1], x,y)e^{\lambda_1 x+ \lambda_2 y}\, \td x \td y\\
&=& \sum_{i=0}^{\IF  }
\int_{\R^2} \pk{\exists_{t \in [i,(i+1)] }: 	\begin{array}{ccc}
		W_1(t) - t > x \\
	W_2(t) - at> y
	\end{array}
}
e^{  \lambda_1 x+ \lambda_2 x_2 } \,  \td x \td y  \\
&\le & \sum_{i=0}^{\IF  }
\int_{\R^2} \pk{\sup_{t \in [i,(i+1)] } 			(W_1(t) - t) > x,
	\sup_{t \in [i,(i+1)] }		(W_{2}(t)- at)> y
}
e^{  \lambda_1 x+ \lambda_2 x_2 } \,  \td x \td y  \\
&=& \frac{1}{{\lambda_1\lambda_2}} \sum_{i=0}^{\IF }
\E{ e^{  \lambda_1 M_i + \lambda_2 M_i^*   }}.
\EQNY
Using further the independence of increments of the Brownian motion,
 the following equality in distribution (abbreviated as $\EQD$) holds  
\BQNY  (M_i,M_i^*) &=&
\Bigl( \sup_{t \in [i,i+1]} (W_1(t) - t),
\sup_{t \in [i,i+1]} (W_{2}  (t) -at) \Bigr)\\
&\EQD & \Bigl( \sup_{t \in [0,1]} (W_1(t) -t) ,
\sup_{t \in [0,1]} (W_{2}(t) - at) \Bigr) + \Bigl( V_1(i) -i, V_2 (i) -  ai  \Bigr)\\
&=:&  (Q_1,Q_2)
+ \Bigl(V_1 (i)-i , V_2(i) -ai  \Bigr),
\EQNY
with $(V_1,V_2)$ an independent copy of $(W_1,W_2)$. By the definition of $\lambda_1$ and $\lambda_2$ we have
$\lambda_1+ \lambda_2 \rho=1$.
Consequently, since for $\tilde V_{2}$ an independent copy of $V_1$
$$ \lambda_1 V_1(i) + \lambda_2 V_2(i)\,\,  {\stackrel{d}{=}}\,\,  (\lambda_1 +\lambda_2 \rho) V_1(i)
+ \lambda_2 \rho^*  \tilde V_2 (i) = V_1(i) + \lambda_2 \rho^* \tilde V_{2} (i)$$
we obtain 
\BQNY   \ln \E{ e^{   \lambda_1 M_i + \lambda_2 M_i^*  }}
- \ln \E{e^{  \lambda _1 Q_1+ \lambda_2 Q_2}} &=&
\frac{i}{2}+ \frac{(a- \rho)^2}{2(1- \rho^2)} i -  \frac{1- a\rho}{1- \rho^2}i
- \frac{a-\rho}{1- \rho^2} ai \\
&=&
-  \frac i {2(1- \rho^2)} [ 2 + 2a^2 - (1- \rho^2)- (a-\rho)^2]\\
&=& - i {\frac{\kappa}{2}} ,
\EQNY
\cE{where $\kappa=\frac{ 1 \kkk{-} 2 a \rho+ a^2}{1- \rho^2} > 1$.}
Consequently, by the monotone convergence theorem
$$  \limit{T } \int_{\R^2} h(T,x,y)e^{\lambda_1 x+ \lambda_2 y}\, \td x \td y =
 \int_{\R^2} \widetilde h(x,y)e^{\lambda_1 x+ \lambda_2 y} \, dx dy \in (0,\IF)$$
implying the claim.\\
$ii)$ Next suppose that $a \le \rho$.
Again  by \nelem{singA} the proof that $\lim_{T\to\infty}I(T)\in (0,\infty)$
follows if we show  that {
$$ \int_{\R}e^{x} \,\pk{ \sup_{t\ge 0 } (W_1(t)- t)> x} \, \td x< \IF.$$
In view of the fact that $\pk{ \sup_{t\ge 0 } (W_1(t)- t)> x}=e^{-2x}$ for $x\geq 0$ we have that
$$\int_{\R}e^{x}\, \pk{ \sup_{t\ge 0 } (W_1(t)- t)> x} \, \td x= 2$$}
implying that for $a< \rho$
\BQNY
 \int_{\R^2} \widetilde h(x,y) e^{x- \frac{y^2- 2 y(c_2 - c_1 \rho)}{2(1-\rho^2)}}\, \td x \td y
 &=& {2}
\int_{\R}  e^{- \frac{y^2- 2 y(c_2 - c_1 \rho)}{2(1-\rho^2)}}\td y\\
& =& 2  \sqrt{2\pi (1- \rho^2)} e^{ \frac{ (c_2- \rho c_1)^2}{2(1- \rho^2)}}
\EQNY
and for $a=\rho$
\BQNY
\lefteqn{
 \int_{{\R\times (-\infty,0)}} \widetilde h(x,y) e^{x- \frac{y^2- 2 y(c_2 - c_1 \rho)}{2(1-\rho^2)}}\, \td x \td y}\\
&= & {2}
\int_{0}^\IF   e^{- \frac{y^2- 2 y(c_1 \rho- c_2)}{2(1-\rho^2)}} \td y\\
&=&2
e^{ \frac{ (c_2- \rho c_1)^2}{2(1- \rho^2)}}
\int_{0}^\IF   e^{- \frac{(y- (c_1 \rho- c_2))^2}{2(1-\rho^2)}} \td y \\
&=&
2 e^{ \frac{ (c_2- \rho c_1)^2}{2(1- \rho^2)}} \sqrt{2 \pi (1- \rho^2)} \Phi((c_1 \rho- c_2)/\sqrt{1- \rho^2}),
\EQNY
where $\Phi$ is the standard normal distribution, hence the proof follows easily.
\QED

\subsection{Proof of Theorem \ref{Th_0}}
The proof is based on the observation that each case i), ii), and iii)
makes reduction of the original problem to a simpler one, which can be solved.

In the light of  \cite{Michna11} and \cite{Michna15} 
for $x\ge0$ we have
$$
\pk{\sup_{t\in [0,T]} (Z(t)-ct)>x}={\mathcal{L}}(c, T, x)\,,
$$
where ${\mathcal{L}}(c, T, x)$ is defined in (\ref{calL}).

The proof of case i) follows now easily since $c_1> c_2$ implies that
$$ \psi_Z(x,y) = \pk{ \sup_{t\in [0,T]} (Z(t)- c_1 t)> x}$$
for any $x\ge y \ge 0$. The case iii) follows with similar arguments, therefore we prove next only the remaining claim.

ii) If $0\leq x< y < x+T(c_1-c_2)$, then
\begin{eqnarray*}
\psi_{Z}(x,y)&=&\pk{\sup_{t\in[0, T]} (Z(t)-c(t)) >y},
\end{eqnarray*}
with (set below $\delta= c_1 -c_2>0$)
\begin{equation}\label{d}
c(t)=
\left\{\begin{array}{ll}
c_2 t &\mbox{if }\, t\in\left[0, \frac{y-x}{\delta}\right]\\
c_1t+x-y&\mbox{if }\, t\in (\frac{y-x}{\delta},T].
\end{array}
\right.
\end{equation}
Hence, following  \cite{Mic17}[Thm 4] we have with $\xi= (y-x)/\delta \in (0,T)$
\begin{eqnarray*}
\psi_{Z}(x,y)&=&
\mathcal{L}\left(c_2,  \xi , y\right)+\int_{0}^\infty \mathcal{L}\left(c_1, T- \xi, z\right)
f(y+c_2-z, T)\td z\\
&&\,\,\,\,\,\,-\int_{0}^\infty z\,\mathcal{L}\left(c_1, T- \xi, z\right)\td z\int_{0}^{\xi}\frac{f(y+c_2 s,s)}
{\xi-s}f\left(c_2\left(\xi-s\right)-z,\xi-s\right)\td s
\end{eqnarray*}
establishing the proof. \QED

\subsection{Proof of Theorem \ref{Th_01}}
The proof of Theorem \ref{Th_01} is analogous to the proof of
Theorem \ref{Th_0}.
\QED

\subsection{Proof of Theorem \ref{Th_3}}
 We focus only on the case $c_1,c_2\ge0$,
since other scenarios follow by similar arguments.
Using that (recall $W_i^*(t)= W_i(t)- c_it$), for any $u>0$ and $x>0$
\[
\pk{u^2(1-\tau_{sim}(u))> x|\tau_{sim}(u)\le 1}
=\frac{\pk{\exists_{t\in [0,1-\frac{x}{u^2}]}: W_1^*(t) > u, W_2^*(t) > au}}
{\pk{\exists_{t\in [0,1]}: W_1^*(t) > u, W_2^*(t) > au}}
\]
in conjunction with Theorem \ref{Th1}, we are left with
finding the asymptotics of
\BQNY
\pk{\exists_{t\in [0,1-\frac{x}{u^2}]}: W_1^*(t) > u, W_2^*(t)> au}
\EQNY
as $u\to\infty$.
By self-similarity of $W_1, W_2$, we have (recall that we denote $\bar x_u = 1- x/u^2$)
\begin{eqnarray*}
\Pi(u,x) &:=&\pk{\exists_{t\in [0, \bar x_u ]}: W_1^*(t) > u, W_2^*(t) > au} \\
&=&
\pk{\exists_{t\in [0,1]}:\sqrt{ \bar x_u } W_1(t)- c_1  t \bar x_u > u,
\sqrt{\bar x_u } W_2(t)- c_2  t\bar x_u > au}.
\end{eqnarray*}
Consequently, for all $u,x$ positive
\begin{eqnarray}
\Pi(u,x)\ge  \psi\left(\frac{u}{\sqrt{\bar x_u }},\frac{au}{\sqrt{\bar x_u }}\right)\label{pi_u}
\end{eqnarray}
and
\begin{eqnarray}
\Pi(u,x)\le  \psi\left(\frac{u}{\sqrt{\bar x_u }}-\frac{c_1x}{u^2},\frac{au}{\sqrt{\bar x_u }}-\frac{c_2x}{u^2}\right)\label{pi_l}.
\end{eqnarray}
Hence the proof follows by a direct application of  Theorem \ref{Th1} to (\ref{pi_u}) and (\ref{pi_l}).
\QED

	\section*{Acknowledgments}
Partial  supported by  SNSF Grant 200021-175752/1 is kindly acknowledged. KD
was partially supported by NCN Grant No 2015/17/B/ST1/01102 (2016-2019).

\bibliographystyle{ieeetr}
\def\polhk#1{\setbox0=\hbox{#1}{\ooalign{\hidewidth
  \lower1.5ex\hbox{`}\hidewidth\crcr\unhbox0}}}


\end{document}